%% file: 9605005.tex
\font\eusm=eusm10                   


\font\eusms=eusm7                       

\font\eusmss=eusm5                      

\input amstex

\documentstyle{amsppt}
  \magnification=1100
  \hsize=6.2truein
  \vsize=9.0truein
  \hoffset 0.1truein
  \parindent=2em

\NoBlackBoxes

\newcount\mycitestyle \mycitestyle=1

\input mymacros

\notasec
  \newtag{\AxiA}
\newsec{\NotElLem}
  \newtag{\Vball}
 \newproc{\Gammas}
 \newproc{\twonorm}
\newsec{\PropCFGF}
 \newproc{\PropCDef}
 \newproc{\chiIPr}
  \newtag{\chiipr}
  \newtag{\chiiprsa}
  \newtag{\AsinGamma}
  \newtag{\APBP}
  \newtag{\ABDsinGamma}
  \newtag{\AEBE}
  \newtag{\EBEEpBEp}
 \newproc{\LFnoPC}
  \newtag{\chin}
\newsec{\FinMultAbSub}
 \newproc{\Compr}
 \newproc{\AsGetQs}
  \newtag{\pandy}
  \newtag{\AGamma}
  \newtag{\AQBQ}
  \newtag{\ADBGamma}
  \newtag{\ADBD}
 \newproc{\Logomega}
  \newtag{\AQtBQt}
 \newproc{\PandY}
 \newproc{\FreeEntropyMI}
 \newproc{\NoFinMult}
  \newtag{\chinpr}
 \newproc{\Interp}
 \newproc{\LFinf}

\newbib{\DykemaZZInterp}
\newbib{\GeZZSimpleAbelian}
\newbib{\RadulescuZZAmalgSubfactors}
\newbib{\PopaZZPropC}
\newbib{\PopaZZCommModComp}
\newbib{\SzarekZZNets}
\newbib{\VoiculescuZZCircSemiCirc}
\newbib{\VoiculescuZZFreeEntropyII}
\newbib{\VoiculescuZZFreeEntropyIII}

\topmatter

  \title Two applications of free entropy \\
          ***** Revised Version *****
  \endtitle

  \author Kenneth J\. Dykema
  \endauthor

  \date 6 May 1996 \enddate

  \affil 
    Department of Mathematics and Computer Science \\
    Odense Universitet, Campusvej 55 \\
    DK-5230 Odense M \\
    Denmark
  \endaffil

  \address
    Department of Mathematics and Computer Science,
    Odense Universitet, Campusvej 55,
    DK-5230 Odense M,
    Denmark
  \endaddress

  \abstract
    Using Voiculsecu's free entropy, it is shown that the free group factors
    $L(F_n)$ lack property~C of Popa and that they lack finite multiplicity
    abelian subalgebras.
  \endabstract




\endtopmatter

\document \TagsOnRight \baselineskip=18pt

\footnote""{Some of the results in this paper appeared as Odense preprint 1996
Nr\. 3}
  
\vskip3ex
\noindent{\bf Introduction.}
\vskip3ex

In~\cite{\VoiculescuZZFreeEntropyII} and~\cite{\VoiculescuZZFreeEntropyIII},
Voiculescu introduced free entropy for $n$--tuples of self--adjoint elements in
a II$_1$--factor, and used it to prove that free group factors, $L(F_n)$, lack
Cartan 
subalgebras~\cite{\VoiculescuZZFreeEntropyIII}.
In~\cite{\PopaZZPropC}, S\. Popa introduced a property for II$_1$--factors,
called property~C.
(See~\cite{\PopaZZCommModComp} for a paper related
to~~\cite{\PopaZZPropC}.)
Like Property~$\Gamma$ of Murray and von Neumann, this is an asymptotic
commutivity property, but it is formally weaker than property~$\Gamma$.
Factors possessing Cartan subalgebras have property~C.
After Voiculescu's striking result, it is a natural
question whether the factors $L(F_n)$ have property~C.
In this note, we show that they do not.

Liming Ge~\cite{\GeZZSimpleAbelian} used Voiculescu's free entropy to show that
the free group factors $L(F_n)$ for $2\le n<\infty$ lack simple (i.e\. of 
multiplicity one) abelian subalgebras.
We say that an abelian subalgebra, $A$, of a II$_1$--factor $\MvN$ with trace
$\tau$ has {\it finite multiplicity} $m$ if
there are $\xi_1,\ldots,\xi_m\in L^2(\MvN,\tau)$ for which
$$ A\xi_1A+\cdots+A\xi_mA \tag{\AxiA} $$
is dense in $L^2(\MvN,\tau)$ and $m$ is least such that this holds, where
juxtaposition in~(\AxiA) indicates the left
and right actions of $\MvN$ on $L^2(\MvN,\tau)$.
We also give a proof, based on Ge's argument, 
that the free group factors
$L(F_n)$ for $2\le n\le\infty$ have no abelian subalgebras of finite
multiplicity.

\medpagebreak
\vskip1ex
\noindent{\bf Acknowledgement.}
\vskip1ex
Thanks to Sorin Popa for suggesting that Theorem~\LFnoPC{} might be true.
I would also like to thank Liming Ge for showing me an early version
of~\cite{\GeZZSimpleAbelian} and Franz Lehner and Roland Speicher for pointing
out an error in my first Lemma~2.2.

\vskip3ex
\noindent{\bf\S\NotElLem.  Notation and elementary lemmas.}
\vskip3ex

See~\cite{\VoiculescuZZFreeEntropyII} and~\cite{\VoiculescuZZFreeEntropyIII}
for many definitions, notations and theorems associated with free entropy.

In this paper:

\vskip1ex
\noindent{\it Euclidean norms.}
$\Mksa$ means the $k^2$--dimensional Euclidean space of self--adjoint
$k\times k$ matrices over the complex numbers with the Euclidean norm
$\nm{T}_e\eqdef\Tr(A^2)^{1/2}$,
where $\Tr$ is the trace on $\Mksa$ for which $\Tr(1)=k$.
On
$$ (\Mksa)^n=\undersetbrace{n\text{ times}}\to{\Mksa\times\cdots\times\Mksa},
$$
the Euclidean norm is
$\nm{(T_1,\ldots,T_n)}_e\eqdef(\sum_1^n\Tr(T_j)^2)^{1/2}$.

\vskip1ex
\noindent{\it $2$--norms.}
In both $\Mksa$ and in II$_1$--factors, for self--adjoint $x$,
$\nm{x}_2=\tau(x^2)^{1/2}$ where $\tau$
is the trace normalized so that $\tau(1)=1$.
Thus on $\Mksa$ we have $\nm{\cdot}_e=k^{1/2}\nm{\cdot}_2$.

\vskip1ex
\noindent{\it Operator norms}
Without a subscript, $\nm{X}$ means the operator norm of $X$.

\vskip1ex
\noindent{\it Euclidean volumes.}
Free entropy is defined in terms of Euclidean volumes of certain subsets of
$(\Mksa)^n$.
Whenever we have a measurable subset, $X$, of Euclidean space, $E\cong R^n$,
we denote by $\volume_n(X)$ the ($n$--dimensional) Lebesgue measure of $X$.
Thus, for example, what is $\lambda(\Gamma_R(x_1,\cdots,x_n\,;\,m,k,\epsilon))$
in~\cite{\VoiculescuZZFreeEntropyII} and~\cite{\VoiculescuZZFreeEntropyIII} is
here written $\volume_{nk^2}(\Gamma_R(x_1,\cdots,x_n\,;\,m,k,\epsilon))$.

\vskip1ex
\noindent{\it Volumes of Euclidean balls.}
$\Vc_n(r)$ is how we denote the $n$--dimensional volume of a ball of radius $r$
in $R^n$, which is known to be
$$ \Vc_n(r)=\frac{\pi^{n/2}r^n}{\Gamma(1+\frac n2)}. \tag{\Vball} $$

\vskip1ex
\noindent{\it Gamma.}
Let us remark that ``$\Gamma_R$'' always refers to the approximating
subsets defined by Voiculescu, while ``$\Gamma$'' without a subscript always
means the Gamma--function.

\proclaim{Lemma \Gammas}
Let $n,m\in\Naturals$.
Then
$$ \frac1{\Gamma(1+\frac n2)\Gamma(1+\frac m2)}\le\frac{2^{(n+m)/2}}
{\Gamma(1+\frac{n+m}2)}. $$
\endproclaim
\demo{Proof}
Taking the crossed product of balls we see
$\Vc_n(2^{-1/2})\Vc_m(2^{-1/2})\le\Vc_{n+m}(1)$.
Use~(\Vball).
\QED

\proclaim{Lemma \twonorm}
Let $0<\beta<1$, $n\in\Naturals$ and $0\le t_1,\ldots,t_N\le\beta$ be such that
$\sum_{i=1}^Nt_i=1$.
Then
$$ \sum_{i=1}^Nt_i^2\le\beta+2\beta^2. $$
\endproclaim
\demo{Proof}
Since $(t_1+c)^2+(t_2-c)^2\ge t_1^2+t_2^2$ if $c\ge0$ and $t_1\ge t_2$, it
follows that the supremum of $\sum_1^Nt_i^2$ over allowable $t_j$
is $[\frac1\beta]\beta^2+(1-\beta[\frac1\beta])^2$.
Hence if $\beta=\frac1K$ for $K\in\Naturals$ then $\sum_1^Nt_i^2\le\frac1K$.
For general $0<\beta<1$ let $K=[\frac1\beta]$.
Then $\frac1{K+1}<\beta\le\frac1K$ so also $0<t_i\le\frac1K$ and thus
$$ \sum_1^Nt_i^2\le\frac1K\le\beta+\frac1{K(K+1)}
\le\beta+\beta^2(1+\frac1K). $$
\QED

\vskip3ex
\noindent{\bf\S\PropCFGF.  Property C and free entropy.}
\vskip3ex
\nopagebreak
\proclaim{Definition \PropCDef}\rm(\cite{\PopaZZPropC}).
A II$_1$--factor $\MvN$ with tracial state $\tau$ is said to have {\it
property C}\/ if $\forall x_1,\ldots,x_n\in\MvN$
$\forall\alpha>0$ $\exists y_1,\ldots,y_m\in\MvN$ such that
\roster
\item"(a)" $\forall1\le\iota\le n$ the distance with respect to
$\nm{\;}_2$ from $x_\iota$
to $\lspan\{y_1,\ldots,y_m\}$ is less than $\alpha$,
\item"(b)" $\forall\beta>0$ $\exists A_1,\ldots,A_m$, mutually commuting
abelian subalgebras of $\MvN$, none of which contains a minimal projection of
trace larger than $\beta$,
and such that for each $j$ the distance with respect to
$\nm{\;}_2$ from $y_j$ to
$A_j'\cap\MvN$ is less than $\beta$.
\endroster
\endproclaim

\proclaim{Lemma \chiIPr}
For $\MvN$ a II$_1$--factor with tracial state $\tau$, for
self--adjoint $x_1,\ldots,x_n\in\MvN$ and for $0<\alpha<1$, suppose
$y_1,\ldots,y_m$ are such that (a) and~(b) of
Definition~\PropCDef{} hold.
Then there is $\beta_0>0$ such that whenever $0<\beta<\beta_0$, whenever
$A_1,\ldots,A_m$
are the mutually commuting abelian subalgebras of $\MvN$ as in~(b) of
Definition~\PropCDef{} and whenever
for each $1\le j\le m$ we take $N_j\in\Naturals$ and
self--adjoint projections $p_{j,1},\ldots,p_{j,N_j}\in A_j$, each of trace
$\le\beta$, such that $\sum_{l=1}^{N_j}p_{j,l}=1$,
then letting $b=\max(\nm{x_1}_2,\ldots,\nm{x_n}_2)+1$ we have
$$ \chi(x_1,\ldots,x_n:y_1,\ldots,y_m,(p_{j,l})
\Sb1\le j\le m\\1\le l\le N_j\endSb)\le C_1+(n-1-6m\beta)\log\alpha,
\tag{\chiipr} $$
where $C_1$ is a constant depending only on $b$ and $n$.
\endproclaim
\demo{Proof}
At the cost of doubling $m$ we may suppose each $y_j$ is self--adjoint.
We will assume each of $y_1,\ldots,y_m$ is self-adjoint and prove under this
additional hypothesis that 
$$ \chi(x_1,\ldots,x_n:y_1,\ldots,y_m,(p_{j,l})
\Sb1\le j\le m\\1\le l\le N_j\endSb)\le C_1+(n-1-3m\beta)\log\alpha,
\tag{\chiiprsa} $$
which will then prove the lemma.
Let $c_{\iota,j}\in\Reals$ be such that
$$ \forall1\le\iota\le n\quad\nm{x_\iota-\sum_{j=1}^mc_{\iota,j}y_j}_2<\alpha.
$$
Take $\beta>0$ and let $A_1,\ldots,A_m$ be the mutually commuting subalgebras
with projections $p_{j,l}$ as in the statement of the lemma.
If $B_j=\Cpx p_{j,1}+\Cpx p_{j,2}+\cdots+\Cpx p_{j,N_j}$ then the distance with
repect to $\nm{\;}_2$ from
$y_j$ to $B_j'\cap\MvN$ is $\nm{y_j-\sum_{l=1}^{N_j}p_{j,l}y_jp_{j,l}}_2$,
which by hypothesis must be $<\beta$.
Hence choosing $\beta_0$ small enough and $\beta<\beta_0$ we get
$$ \forall1\le\iota\le n\qquad
\nm{x_\iota-\sum_{j=1}^mc_{\iota,j}\sum_{l=1}^{N_j}p_{j,l}y_jp_{j,l}}_2
<\alpha+\beta\sum_{j=1}^m|c_{\iota,j}|<2\alpha. $$
Let $q_1,\ldots,q_N$ be the minimal projections of the algebra generated by
$\{p_{j,l}\mid1\le j\le m,\,1\le l\le N_j\}$.
Let $I(j,l)\subseteq\{1,\ldots,N\}$ be such that
$$ p_{j,l}=\sum_{s\in I(j,l)}q_s. $$
By an easy case of~\scite{\VoiculescuZZFreeEntropyIII}{1.8}, we have that
$$ \chi(x_1,\ldots,x_n:y_1,\ldots,y_m,
(p_{j,l})\Sb1\le j\le m\\1\le l\le N_j\endSb)
=\chi(x_1,\ldots,x_n:y_1,\ldots,y_m,q_1,\ldots,q_N) $$
and we estimate the latter quantity.

\proclaim{Claim \chiIPr a}\rm
$\forall R\ge0$
$\exists m',k_0\in\Naturals$ $\exists\epsilon>0$ such that
if
$$
(A_1,\ldots,A_n)\in\Gamma_R(x_1,\ldots,x_n:y_1,\ldots,y_m,q_1,\ldots,q_N;
m',k,\epsilon) \tag{\AsinGamma} $$
then there are $B_1,\ldots,B_m\in\Mksa$ and projections
$Q_1,\ldots,Q_N\in\Mksa$ such that $\sum_1^NQ_s\le1$, $\rank(Q_s)=[\tau(q_s)k]$
and such that, letting
$$ P_{j,l}=\sum_{s\in I(j,l)}Q_s, $$
we have
$$ \forall1\le\iota\le n\quad\nm{A_\iota-\sum_{j=1}^mc_{\iota,j}
\sum_{l=1}^{N_j}P_{j,l}B_jP_{j,l}}_2<3\alpha. \tag{\APBP} $$
\endproclaim
\demo{Proof}
Indeed, if $A_\iota$ are as in~(\AsinGamma), let
$B_1,\ldots,B_m,D_1,\ldots,D_N\in\Mksa$ be such that
$$ (A_1,\ldots,A_n,B_1,\ldots,B_m,D_1,\ldots,D_N)
\in\Gamma_R(x_1,\ldots,x_n,y_1,\ldots,y_m,q_1,\ldots,q_N;m',k,\epsilon).
\tag{\ABDsinGamma} $$
Then by~(\AsinGamma), letting $E_{j,l}=\sum_{s\in I(j,l)}D_s$ and choosing
$m'\ge3$ and $\epsilon$ small enough, we may insist that
$$ \forall1\le\iota\le n\quad\nm{A_\iota-\sum_{j=1}^mc_{\iota,j}
\sum_{l=1}^{N_j}E_{j,l}B_jE_{j,l}}_2<5\alpha/2 \tag{\AEBE} $$
Now there is $\delta>0$ such that if $D_1',\ldots,D_N'\in\Mksa$ satisfy
$\nm{D_s'}\le R$ and
$\nm{D_s'-D_s}_2<\delta$ then letting
$E'_{j,l}=\sum_{s\in I(j,l)}D'_s$ it follows that, for any $\nm{B_j}\le R$,
$$ \forall1\le\iota\le n\quad\nm{\sum_{j=1}^mc_{\iota,j}\sum_{l=1}^{N_j}
(E_{j,l}B_jE_{j,l}-E'_{j,l}B_jE'_{j,l})}_2<\alpha/2. \tag{\EBEEpBEp} $$
It follows from~\scite{\VoiculescuZZFreeEntropyII}{4.3} that if $m'$ is large
enough, if $\epsilon$ is small enough and if $k\ge k_0$ for $k_0$ large enough,
(independently of the choice of $D_1,\ldots,D_N$ such that~(\ABDsinGamma)
holds),
then there are $Q_1,\ldots,Q_N\in\Mksa$, projections, such that
$\sum_1^NQ_s\le1$ and $\forall s$ $\rank(Q_s)=[\tau(q_s)k]$ and
$\nm{Q_s-D_s}_2<\delta$.
Hence, from~(\EBEEpBEp) and~(\AEBE) we have that~(\APBP) holds.
Thus, Claim~\chiIPr a is proved.
\enddemo

Continuing with the proof of Lemma~\chiIPr, let $\Qt_1,\ldots,\Qt_N\in\Mksa$ be
any fixed projections such that $\sum_{s=1}^N\Qt_s\le1$ and
$\rank(\Qt_s)=[\tau(q_s)k]$.
Let
$$ \Pt_{j,l}=\sum_{s\in I(j,l)}\Qt_s. $$
Then by~(\APBP), there is a $k\times k$ unitary $U$, (letting $U$ be such that
$U^*\Qt_sU=Q_s$), for which we have
$$ \forall1\le\iota\le n\quad\nm{UA_\iota U^*-
\sum_{j=1}^mc_{\iota,j}\sum_{l=1}^{N_j}\Pt_{j,l}(UB_jU^*)\Pt_{j,l}}_2<3\alpha.
$$
Consider the linear transformation
$$ F:(\Mksa)^m\to(\Mksa)^n $$
given by
$$ F(X_1,\ldots,X_m)=
\biggl(\sum_{j=1}^mc_{\iota,j}\sum_{l=1}^{N_j}\Pt_{j,l}X_j\Pt_{j,l}\biggr)
_{\iota=1}^n. $$
Thus every $(A_1,\ldots,A_n)$ as in~(\AsinGamma) is within Euclidean distance
$3\alpha n^{1/2}k^{1/2}$ from a unitary conjugate of the range
space of $F$.
If $U,U_1\in\Uc_k$ are such that $\nm{U-U_1}\le\gamma$ then since
$\nm{A_\iota}_2\le b-1+\epsilon\le b$ we have
$$ \align
\nm{&(U^*A_1U,\ldots,U^*A_nU)-(U^*_1A_1U_1,\ldots,U^*_1A_nU_1)}_e\le \\
&\le k^{1/2}n^{1/2}\max_{1\le\iota\le n}\nm{U^*A_\iota U-U_1^*A_\iota U_1}_2
\le2bk^{1/2}n^{1/2}\gamma. \endalign $$
Hence if $(U_\lambda)_{\lambda\in\Lambda}$ is a $\gamma$--net for $\Uc_k$ with
respect to the metric arising from the operator norm,
then $(A_1,\ldots,A_n)$ is within Euclidean
distance $n^{1/2}k^{1/2}(3\alpha+2b\gamma)$ of
$$ \bigcup_{\lambda\in\Lambda}(U_\lambda,\ldots,U_\lambda)\text{Range}(F)
(U_\lambda^*,\ldots,U_\lambda^*). $$
By results of Szarek~\cite{\SzarekZZNets}, there is a $\gamma$--net in $\Uc_k$
with cardinality $|\Lambda|\le(\frac C\gamma)^{k^2}$, where $C$ is a universal
constant.
In the following computation, let $\Gamma_R$ denote the right--hand--side
of~(\AsinGamma).
Since the Euclidean norm of $(A_1,\ldots,A_n)$ is no greater than
$k^{1/2}n^{1/2}b$, letting $\Bc$ be the Euclidean ball in Range$(F)$ of
radius
$k^{1/2}n^{1/2}(3\alpha+b(2\gamma+1))$,
it follows that every point in $\Gamma_R$
is connected to
$(U_\lambda,\ldots,U_\lambda)\Bc(U_\lambda^*,\ldots,U_\lambda^*)$, for some
$\lambda\in\Lambda$ by a line segment of Euclidean length
$\le k^{1/2}n^{1/2}(3\alpha+2b\gamma)$ which is normal to
$(U_\lambda,\ldots,U_\lambda)\Bc(U_\lambda^*,\ldots,U_\lambda^*)$.
Let $d_k=\rank(F)$.
Then
$$ \align
d_k&\le\sum_{j=1}^m\sum_{l=1}^{N_j}\rank(\Pt_{j,l})^2\\
&\le k^2\sum_{j=1}^m\sum_{l=1}^{N_j}\tau(p_{j,l})^2.
\endalign $$
Hence by Lemma~\twonorm, $d_k\le 3m\beta k^2$.
Then
$$ \volume_{nk^2}(\Gamma_R)\le
|\Lambda|\cdot\Vc_{d_k}(k^{1/2}n^{1/2}(3\alpha+b(2\gamma+1)))\cdot
\Vc_{nk^2-d_k}(k^{1/2}n^{1/2}(3\alpha+2b\gamma)).
$$
Letting $\gamma=\alpha/b$ we have
$$ \align
\volume_{nk^2}(\Gamma_R)&\le
\left(\frac{Cb}\alpha\right)^{k^2}
\cdot\frac{(\pi kn)^{d_k/2}(5\alpha+b)^{d_k}}{\Gamma(1+\frac{d_k}2)}
\cdot\frac{(\pi kn)^{(nk^2-d_k)/2}(5\alpha)^{nk^2-d_k}}
{\Gamma(1+\frac{nk^2-d_k}2)} \\
&\le\frac{(Cb)^{k^2}(\pi kn)^{nk^2/2}(6b)^{d_k}5^{nk^2-d_k}\alpha^{(n-1)k^2-d_k}
2^{nk^2/2}}{\Gamma(1+\frac{nk^2}2)} \\
&\le\frac{((6b)^n5^nCb)^{k^2}(2\pi kn)^{nk^2/2}\alpha^{(n-1)k^2-d_k}}
{\Gamma(1+\frac{nk^2}2)},
\endalign $$
where we used successively~(\Vball) and Szarek's result about $|\Lambda|$;
Lemma~\Gammas; the fact that $0\le d_k\le nk^2$.
Thus using Stirling's formula we get
$$ \align
\frac n2&\log k+k^{-2}\log\volume_{nk^2}(\Gamma_R)\le \\
&\le n\log(30b)+\log(Cb)+\frac n2\log(2\pi n)
+(n-1-(d_k/k^2))\log\alpha+c_k, \endalign $$
where $\lim_{k\to\infty}c_k=0$.
Using that $d_k\le 3m\beta k^2$ and letting $k\to\infty$ gives~(\chiiprsa).
\QED

\proclaim{Theorem \LFnoPC}
For each $n\in\{2,3,4,\ldots,\infty\}$, $L(F_n)$ does not have property C.
\endproclaim
\demo{Proof}
Suppose for contradiction that $L(F_n)$ has property C.
We have by~\scite{\VoiculescuZZFreeEntropyII}{4.5}
and~\scite{\VoiculescuZZFreeEntropyII}{5.4} that for a free family of
semicircular elements, $x_1,\ldots,x_{n'}$,
$$ \chi(x_1,\ldots,x_{n'})>-\infty. \tag{\chin} $$
If $n\in\Naturals$, $n\ge2$, then letting $x_1,\ldots,x_n$ be a free family of
semicircular elements generating $L(F_n)$, using~(\chiiprsa) with $\beta$ and
$\alpha$ sufficiently small, from~(\chiipr) and~(\chin) we obtain
$$ \chi(x_1,\ldots,x_n:z_1,\ldots,z_K)<\chi(x_1,\ldots,x_n) $$
for some $z_1,\ldots,z_K\in L(F_n)$, which is a contradiction
to~\scite{\VoiculescuZZFreeEntropyIII}{1.8}.

If $n=\infty$, let $(x_\iota)_{\iota=1}^\infty$ be a free family of
semicircular elements generating $L(F_\infty)$ and
take any $n'\in\Naturals$, $n'\ge2$.
Assuming property C holds, in like manner to above we obtain
$$ \chi(x_1,\ldots,x_{n'}:z_1,\ldots,z_K)<\chi(x_1,\ldots,x_{n'}) $$
for some $z_1,\ldots,z_K\in L(F_\infty)$.
Now letting $E_l$ denote the conditional expectation from $L(F_\infty)$ to
$\{x_1,\ldots,x_l\}$, by~\scite{\VoiculescuZZFreeEntropyIII}{1.5} we have
for some $l\in\Naturals$, $l\ge n'$, that
$$ \chi(x_1,\ldots,x_{n'}:E_l(z_1),\ldots,E_l(z_K))<\chi(x_1,\ldots,x_{n'}). $$
Hence by~\scite{\VoiculescuZZFreeEntropyIII}{1.8}
and~\scite{\VoiculescuZZFreeEntropyII}{2.3} we have
$$ \align
\chi(x_1,\ldots,x_l)&=\chi(x_1,\ldots,x_l:z_1,\ldots,z_K) \\
&\le\chi(x_1,\ldots,x_{n'}:z_1,\ldots,z_K)+\chi(x_{n'+1},\ldots,x_l)\\
&<\chi(x_1,\ldots,x_{n'})+\chi(x_{n'+1})+\cdots+\chi(x_l)\\
&\le\chi(x_1)+\cdots+\chi(x_l), \endalign $$
which contradicts~\scite{\VoiculescuZZFreeEntropyII}{5.4}.
\QED

\vskip3ex
\noindent{\bf\S\FinMultAbSub. Finite multiplicity abelian subalgebras and free
entropy.}
\vskip3ex

\proclaim{Proposition \Compr}
Suppose a II$_1$--factor, $\MvN$, has an abelian subalgebra, $A$, of finite
multiplicity $\le m$.
Then for every $t\in\Real_{>0}$, $\MvN_t$ has an abelian subalgebra of
multiplicity $\le m$.
\endproclaim
\demo{Proof}
It will suffice to show for all $0<t<1$ and for all $t\in\Naturals$
that if $A$ has multiplicity $\le m$ then $\MvN_t$ has an abelian subalgebra of
finite multiplicity $\le m$.
Let $\xi_1,\ldots,\xi_m\in L^2(\MvN,\tau)$ be such that
$A\xi_1A+\cdots+A\xi_mA$ is dense in $L^2(\MvN,\tau)$.
For $0<t<1$, let $p\in A$ be a projection for which $\tau(p)=t$.
Then
$$ \align
L^2(p\MvN p,\tau(p)^{-1}\restrict_{p\MvN p})&=p L^2(\MvN,\tau)p \\
&=p\overline{(A\xi_1A+\cdots A\xi_m A)}p \\
&=\overline{pA(p\xi_1p)Ap+\cdots pA(p\xi_mp)Ap} \endalign $$
so $pA$ has multiplicity $\le m$ in $p\MvN p$.
For $t=n\in\Naturals$, $t\ge2$ we have $\MvN_t=\MvN\otimes M_n(\Cpx)$.
Let $(e_{ij})_{1\le i,j\le n}$ be a system of matrix units for $M_n(\Cpx)$.
Then $L^2(M_n(\Cpx),\tau_n)$ has orthonormal basis
$\{n^{1/2}\eh_{ij}\mid1\le i,j\le n\}$.
Let $D=\lspan\{e_{ii}\mid 1\le i\le n\}\subseteq M_n(\Cpx)$.
Consider the abelian subalgebra, $A\otimes D\subseteq\MvN_t$.
Let $\xi_k'=\sum_{i,j=1}^n\xi_k\otimes\eh_{ij}$, so that
$$ (a_1\otimes e_{ii})\xi'_k(a_2\otimes e_{jj})=a_1\xi_ka_2\otimes\eh_{ij} $$
and thus
$$ (A\otimes D)\xi_1'(A\otimes D)+\cdots+(A\otimes D)\xi_m'(A\otimes D) $$
is dense in $L^2(\MvN,\tau)\otimes L^2(M_n(\Cpx),\tau_n)$.
Hence $A\otimes D$ has multiplicity $\le m$ in $\MvN\otimes M_n(\Cpx)$.
\QED

\proclaim{Lemma \AsGetQs}
Let $\MvN$ be a II$_1$--factor with normalized trace $\tau$ and
take $x_1,\ldots,x_n\in\MvNsa$.
Suppose $\omega>0$ and there are $y_1,\ldots,y_K\in\MvNsa$,
projections $p_1,\ldots,p_N\in\MvNsa$ whose sum is $1$ and scalars
$\lambda_{r,s}^{(\iota,l)}\in\Cpx$
(for $1\le\iota\le n$,$\,1\le l\le K$,$\,1\le r,s\le N$) for which
$\lambda_{r,s}^{(\iota,l)}=\overline{\lambda_{s,r}^{(\iota,l)}}$ such that
$$ \forall1\le\iota\le n\qquad\nm{x_\iota-\sum_{l=1}^K\sum_{r,s=1}^N
\lambda_{r,s}^{(\iota,l)}p_ry_lp_s}_2<\omega. \tag{\pandy} $$
Then $\forall R>1$ there are $k_1,m'\in\Naturals$ and $\epsilon>0$ such that
if $k\ge k_1$ and if
$$ (A_1,\ldots,A_n)\in
\Gamma_R
(x_1,\ldots,x_n:p_1,\ldots,p_N,y_1,\ldots,y_K;m',k,\epsilon) \tag{\AGamma} $$
then there are orthogonal, self--adjoint projections
$Q_1,\ldots,Q_N\in\Mksa$, with $\rank(Q_j)=[\tau(p_j)k]$ and there are
$B_1,\ldots,B_K\in\Mksa$ such that
$$ \forall1\le\iota\le n\qquad\nm{A_\iota-\sum_{l=1}^K\sum_{r,s=1}^N
\lambda_{r,s}^{(\iota,l)}Q_rB_lQ_s}_2<3\omega. \tag{\AQBQ} $$
\endproclaim
\demo{Proof}
In light of~(\pandy), taking $m_0=3$ there is $\epsilon_0>0$ such that if
$k\in\Naturals$ and if
$$ (A_1,\ldots,A_n,D_1,\ldots,D_N,B_1,\ldots,B_K)\in
\Gamma_R(x_1,\ldots,x_n,p_1,\ldots,p_N,y_1,\ldots,y_K;m_0,k,\epsilon_0)
\tag{\ADBGamma} $$
then
$$ \forall1\le\iota\le n\qquad\nm{A_\iota-\sum_{l=1}^K\sum_{r,s=1}^N
\lambda_{r,s}^{(\iota,l)}D_rB_lD_r}_2<2\omega. \tag{\ADBD} $$
Then there is $\delta>0$ such that if $D_1',\ldots,D_N'\in\Mksa$,
$\nm{D_r'}\le R$ satisfy
$\forall1\le j\le N$ $\nm{D_j-D_j'}_2<\delta$ then for any $\nm{B_l}\le R$,
$$ \forall1\le\iota\le n\qquad\nm{\sum_{l=1}^K\sum_{r,s=1}^N
\lambda_{r,s}^{(\iota,l)}(D_rB_lD_s-D_r'B_lD_s')}_2<\omega. $$
It follows from~\scite{\VoiculescuZZFreeEntropyII}{4.3} that there are
$\epsilon_1>0$, $m_1,k_1\in\Naturals$ such that if $k\ge k_1$ and if
$$ (D_1,\ldots,D_N)\in\Gamma_R(p_1,\ldots,p_N;m_1,k,\epsilon_1) $$
then there are orthogonal, self--adjoint projections $Q_1,\ldots,Q_N\in\Mksa$,
$\rank(Q_j)=[\tau(p_j)k]$, such that $\forall 1\le j\le N$
$\nm{D_j-Q_j}_2<\delta$.

Putting these facts together, we can prove the lemma.
Indeed, let $\epsilon=\min(\epsilon_0,\epsilon_1)$ and $m'=\max(m_0,m_1)$.
If $A_\iota$ are as in~(\AGamma) then there are $D_j$ and $B_l$ such
that~(\ADBGamma) holds, also with $m'$ and $\epsilon$ replacing $m_0$ and
$\epsilon_0$.
Then~(\ADBD) holds and we can find $Q_j$ as desired so that~(\AQBQ) holds.
\QED

\proclaim{Lemma \Logomega}
Suppose $x_1,\ldots,x_n$, $\omega>0$, $y_l$, $p_j$ and
$\lambda_{r,s}^{(\iota,l)}$
satisfy the hypotheses of Lemma~\AsGetQs{} and let
$b=\max(\nm{x_1}_2,\ldots,\nm{x_n}_2)+1$.
Let $m'\in\Naturals$ and $\epsilon>0$ be as found in Lemma~\AsGetQs{},
assume without loss of generality that $\epsilon\le1$ and let
$$ R>\max(\nm{x_1},\ldots,\nm{x_n},\nm{y_1},\ldots,\nm{y_K},1). $$
If also $\omega\le b$ and $\omega<1$ then
$$ \chi_R(x_1,\ldots,x_n:y_1,\ldots,y_K,p_1,\ldots,p_N;m',\epsilon)
\le C_1+(n-K-1)\log\omega, $$
where $C_1\in\Reals$ depends only on $n$ and $b$.
\endproclaim
\demo{Proof}
Let $k_1$ be as found in Lemma~\AsGetQs{} and let $k\ge k_1$.
Let $\Qt_1,\ldots,\Qt_N$ be any fixed, orthogonal projections in $\Mksa$, with
$\rank(\Qt_j)=[\tau(p_j)k]$.
By Lemma~\AsGetQs, for every $(A_1,\ldots,A_n)$ as in~(\AGamma) there are
$B_1,\ldots,B_K\in\Mksa$, each of norm $\le R$, and there is
$U\in\Uc_k$ (the group of $k\times k$ unitaries) such that
$$ \forall1\le\iota\le n\qquad\nm{A_\iota-\sum_{l=1}^K\sum_{r,s=1}^N
\lambda_{r,s}^{(\iota,l)}U\Qt_rU^*B_lU\Qt_sU^*}_2<3\omega. \tag{\AQtBQt} $$
Consider the linear transformation
$$ F:(\Mksa)^K\rightarrow(\Mksa)^n $$
given by
$$ F(C_1,\ldots,C_K)=\bigg(\sum_{l=1}^K\sum_{r,s=1}^N
\lambda_{r,s}^{(\iota,l)}\Qt_rC_l\Qt_s\biggr)_{\iota=1}^n. $$
Then~(\AQtBQt) implies that there is $U\in\Uc_k$ such that the point
$(U^*A_1U,\ldots,U^*A_nU)$ is Euclidean distance at most
$3\omega n^{1/2}k^{1/2}$ away from the range space of $F$.
If $U,U_1\in\Uc_k$ are such that $\nm{U-U_1}\le\gamma$ then since
$\nm{A_\iota}_2\le b-1+\epsilon\le b$ we have
$$ \align
\nm{&(U^*A_1U,\ldots,U^*A_nU)-(U^*_1A_1U_1,\ldots,U^*_1A_nU_1)}_e\le \\
&\le k^{1/2}n^{1/2}\max_{1\le\iota\le n}\nm{U^*A_\iota U-U_1^*A_\iota U_1}_2
\le2bk^{1/2}n^{1/2}\gamma. \endalign $$
Hence if $(U_\lambda)_{\lambda\in\Lambda}$ is a $\gamma$--net for $\Uc_k$ with
respect to the metric arising from the operator norm,
then $(A_1,\ldots,A_n)$ is within Euclidean
distance $n^{1/2}k^{1/2}(3\omega+2b\gamma)$ of
$$ \bigcup_{\lambda\in\Lambda}(U_\lambda,\ldots,U_\lambda)\text{Range}(F)
(U_\lambda^*,\ldots,U_\lambda^*). $$
By results of Szarek~\cite{\SzarekZZNets}, there is a $\gamma$--net in $\Uc_k$
with cardinality $|\Lambda|\le(\frac C\gamma)^{k^2}$, where $C$ is a universal
constant.
In the following computation, let $\Gamma_R$ denote the right--hand--side
of~(\AGamma).
Since the Euclidean norm of $(A_1,\ldots,A_n)$ is no greater than
$k^{1/2}n^{1/2}b$, letting $\Bc$ be the Euclidean ball in Range$(F)$ of
radius
$k^{1/2}n^{1/2}(3\omega+b(2\gamma+1))$,
it follows that every point in $\Gamma_R$
is connected to
$(U_\lambda,\ldots,U_\lambda)\Bc(U_\lambda^*,\ldots,U_\lambda^*)$, for some
$\lambda\in\Lambda$ by a line segment of Euclidean length
$\le k^{1/2}n^{1/2}(3\omega+2b\gamma)$ which is normal to
$(U_\lambda,\ldots,U_\lambda)\Bc(U_\lambda^*,\ldots,U_\lambda^*)$.
Let $d_k=\rank(F)$, so $d_k\le nk^2$ and $d_k\le mk^2$.
Then
$$ \volume_{nk^2}(\Gamma_R)\le
|\Lambda|\cdot\Vc_{d_k}(k^{1/2}n^{1/2}(3\omega+b(2\gamma+1))\cdot
\Vc_{nk^2-d_k}(k^{1/2}n^{1/2}(3\omega+2b\gamma)). $$
Set $\gamma=\omega/b$, giving
$$ \align
\volume_{nk^2}(\Gamma_R)&\le
\left(\frac{Cb}\omega\right)^{k^2}
\cdot\frac{(\pi kn)^{d_k/2}(5\omega+b)^{d_k}}{\Gamma(1+\frac{d_k}2)}
\cdot\frac{(\pi kn)^{(nk^2-d_k)/2}(5\omega)^{nk^2-d_k}}
{\Gamma(1+\frac{nk^2-d_k}2)} \\
&\le\frac{(Cb)^{k^2}(\pi kn)^{nk^2/2}(6b)^{d_k}5^{nk^2-d_k}\omega^{(n-1)k^2-d_k}
2^{nk^2/2}}{\Gamma(1+\frac{nk^2}2)} \\
&\le\frac{((6b)^n5^nCb)^{k^2}(2\pi kn)^{nk^2/2}\omega^{(n-K-1)k^2}}
{\Gamma(1+\frac{nk^2}2)},
\endalign $$
where we used successively~(\Vball) and Szarek's result about $|\Lambda|$;
Lemma~\Gammas; the facts that $(n-1)k^2-d_k\ge(n-K-1)k^2$ and $\omega<1$.
Thus using Stirling's formula we get
$$ \align
\frac n2&\log k+k^{-2}\log\volume_{nk^2}(\Gamma_R)\le \\
&\le n\log(30b)+\log(Cb)+\frac n2\log(2\pi n)-\frac12\log(2\pi)
+(n-K-1)\log\omega+c_k, \endalign $$
where $\lim_{k\to\infty}c_k=0$.
Taking the limit as $k\to\infty$ gives
$$\chi_R(x_1,\ldots,x_n:p_1,\ldots,p_N,y_1,\ldots,y_K;m',\epsilon)
\le\,C_1+(n-K-1)\log\omega $$
as required.
\QED

\proclaim{Lemma \PandY}
Suppose $\MvN$ is a II$_1$--factor with abelian subalgebra $A$ having finite
multiplicity $\le m$.
Let $x_1,\ldots,x_n\in\MvNsa$ and $\omega>0$.
Let $K=2m$.
Then there is $N\in\Naturals$ such that the hypotheses of Lemma~\AsGetQs{} are satisfied, namely,
there are $y_1,\ldots,y_{2m}$, $p_1,\ldots,p_N$ and $\lambda_{r,s}^{(\iota,l)}$
such that~(\pandy) holds.
\endproclaim
\demo{Proof}
Recall that an element of $\Hil=L^2(\MvN,\tau)$ is said to be self--adjoint if
it is fixed by the involution $J$ defined by $J\xh=(x^*)\hat{\,}$, where
$x\mapsto\xh$ is the defining mapping $M\to L^2(M,\tau)$.
There are self--adjoint $\xi_1,\ldots,\xi_{2m}\in\Hil$
such that $A\xi_1A+\cdots+A\xi_{2m}A$ is dense in $\Hil$.
Hence there is $M\in\Naturals$ and there are
$a_1^{(\iota,l)},\ldots,a_M^{(\iota,l)},
b_1^{(\iota,l)},\ldots,b_M^{(\iota,l)}\in A$ 
(for $1\le l\le 2m$, $\,1\le\iota\le n$) such that
$$ \forall1\le\iota\le n\qquad\nm{\xh_\iota-\sum_{l=1}^{2m}\sum_{k=1}^M
a_k^{(\iota,l)}\xi_lb_k^{(\iota,l)}-\xh_\iota}_\Hil<\omega/3. $$
Since $(\MvN_{\text{\it s.a\.}})\hat{\,}$ is dense in
$\Hil_{\text{\it s.a\.}}$,
there are $y_1,\ldots,y_{2m}\in\MvN_{\text{\it s.a\.}}$ such that
$$ \forall1\le\iota\le n\qquad\nm{x_\iota-\sum_{l=1}^{2m}\sum_{k=1}^M
a_k^{(\iota,l)}y_lb_k^{(\iota,l)}}_2<2\omega/3. $$
Since $A$ is abelian there are $N\in\Naturals$ and self--adjoint projections
$p_1,\ldots,p_N\in A$ whose sum is $1$
such that each $a_k^{(\iota,l)}$ and $b_k^{(\iota,l)}$ is approximated
sufficiently well in norm by linear combinations of the $p_j$ so that for some
scalars $\lambda_{r,s}^{(\iota,l)}\in\Cpx$ we have
$$ \forall1\le\iota\le n\qquad\nm{x_\iota-\sum_{l=1}^{2m}\sum_{r,s=1}^N
\lambda_{r,s}^{(\iota,l)}p_ry_lp_s}_2<\omega. $$
Since $x_\iota$ is self--adjoint, we may without penalty take
$\lambda_{r,s}^{(\iota,l)}=\overline{\lambda_{r,s}^{(\iota,l)}}$.
\QED

\proclaim{Theorem \FreeEntropyMI}
Suppose $\MvN$ is a II$_1$--factor with an abelian subalgebra, $A$, having
finite multiplicity $\le m$.
Suppose $x_1,\ldots,x_n$ are self--adjoint elements of $\MvN$ which taken
together generate $\MvN$.
If $n>2m+1$ then
$$ \chi(x_1,\ldots,x_n)=-\infty. $$
\endproclaim
\demo{Proof}
Let $\omega>0$, let $N$, $p_1,\ldots,p_N$ and
$y_1,\ldots,y_{2m}$ be as obtained from Lemma~\PandY{}.
and let $b$ and $R$ be as in Lemma~\Logomega.
Since $\MvN=\{x_1,\ldots,x_n\}''$, by~\scite{\VoiculescuZZFreeEntropyIII}{1.8}
we have
$$ \chi(x_1,\ldots,x_n)=\chi_R(x_1,\ldots,x_n)=
\chi_R(x_1,\ldots,x_n:p_1,\ldots,p_N,y_1,\ldots,y_{2m}), $$
and we will estimate the latter quantity.
Let $\epsilon$ and $m'$
be obtained from Lemma~\AsGetQs, taking without loss of generality
$\epsilon\le1$.
Then for $\omega\le b$, by Lemma~\Logomega{},
$$ \chi(x_1,\ldots,x_n)
\le\chi_R(x_1,\ldots,x_n:p_1,\ldots,p_N,y_1,\ldots,y_{2m};m',\epsilon)
\le C_1+(n-2m-1)\log\omega. $$
Now let $\omega\to0$.
\QED

\proclaim{Corollary \NoFinMult}
For no $n\in\Naturals$, $n\ge2$ does the free group factor $L(F_n)$ contain an
abelian subalgebra of finite multiplicity.
\endproclaim
\demo{Proof}
For $n'\in\Naturals$ and $n'\ge2$,
$L(F_{n'})$ is generated by a free family semicircular elements,
$x_1,\ldots,x_{n'}$
and by~\scite{\VoiculescuZZFreeEntropyII}{4.5}
and~\scite{\VoiculescuZZFreeEntropyII}{5.4} it is seen that
$$ \chi(x_1,\ldots,x_{n'})>-\infty. \tag{\chinpr} $$
Suppose for contradiction that $L(F_n)$ has an abelian subalgebra of
finite multiplicity $m$.
By~\cite{\VoiculescuZZCircSemiCirc}, for each $k\in\Naturals$
$L(F_n)_{1/k}\cong L(F_{1+k^2(n-1)})$, so by Proposition~\Compr{} and for $k$ large
enough we can find $n'\in\Naturals$ such that $n'>2m-1$ and $L(F_{n'})$
has an abelian subalgebra of multiplicity $\le m$.
Theorem~\FreeEntropyMI{} and~(\chinpr) now give a contradiction.
\QED

\proclaim{Remark \Interp}\rm
Clearly, 
(see~\cite{\RadulescuZZAmalgSubfactors},~\cite{\DykemaZZInterp})
the same argument shows that for no $1<t<\infty$ does
the interpolated free group factor
$L(F_t)$ contain an abelian subalgebra of finite multiplicity.
\endproclaim

We now use the method of~\scite{\VoiculescuZZFreeEntropyII}{5.3} to consider
also the free group factor on (countably) infinitely many generators.

\proclaim{Theorem \LFinf}
$L(F_\infty)$ has no abelian subalgebras of finite multiplicity.
\endproclaim
\demo{Proof}
Suppose for contradiction that $L(F_\infty)$ contains an abelian subalgebra of
finite multiplicity $m$.
Let $x_1,x_2,\ldots$ be a free family of semicircular elements generating
$L(F_\infty)$ such that $\nm{x_j}_2=1$.
Let $K=2m$, $n=2m+2$, $b=2$ and let $C_1$ be as in Lemma~\Logomega{} for these
values of $n$ and $b$.
Since $\chi(x_1,\ldots,x_n)>-\infty$, there is $\omega>0$
such that
$C_1+\log\omega<\chi(x_1,\ldots,x_n)$.
Let $y_1,\ldots,y_{2m}$, $p_1,\ldots,p_N$ and $\lambda_{r,s}^{(\iota,l)}$ be as
obtained from Lemma~\PandY.
By Lemma~\Logomega,
$$ \chi(x_1,\ldots,x_n:y_1,\ldots,y_{2m},p_1,\ldots,p_N)<\chi(x_1,\ldots,x_n).
$$
\line{For $n'\in\Naturals$ let $E_{n'}$ be the trace--preserving conditional
expectation from $L(F_\infty)$ onto}
$\{x_1,\ldots,x_{n'}\}''$.
Then by~\scite{\VoiculescuZZFreeEntropyIII}{1.5}, for some $n'\in\Naturals$ we
have
$$ \chi(x_1,\ldots,x_n:E_{n'}(y_1),\ldots,E_{n'}(y_{2m}),E_{n'}(p_1),\ldots,
E_{n'}(p_N))<\chi(x_1,\ldots,x_n) $$
and we may assume without loss of generality that $n'\ge n$.
Thus, by~\scite{\VoiculescuZZFreeEntropyIII}{1.8}
and~\scite{\VoiculescuZZFreeEntropyII}{2.3} we have
$$ \align
\chi(x_1,\ldots,x_{n'})
&=\chi(x_1,\ldots,x_{n'}:E_{n'}(y_1),\ldots,E_{n'}(y_{2m}),E_{n'}(p_1),\ldots,
E_{n'}(p_N)) \\
&\le\chi(x_1,\ldots,x_n:E_{n'}(y_1),\ldots,E_{n'}(y_{2m}),E_{n'}(p_1),\ldots,
E_{n'}(p_N)) \\
&\qquad+\chi(x_{n+1},\ldots,x_{n'}) \\
&<\chi(x_1,\ldots,x_n)+\chi(x_{n+1})+\cdots+\chi(x_{n'}) \\
&\le\chi(x_1)+\cdots+\chi(x_{n'}), \endalign $$
which contradicts~\scite{\VoiculescuZZFreeEntropyII}{5.4}.
\QED

\enddocument

%% file: mymacros.tex

\newcount\theTime
\newcount\theHour
\newcount\theMinute
\newcount\theMinuteTens
\newcount\theScratch
\theTime=\number\time
\theHour=\theTime
\divide\theHour by 60
\theScratch=\theHour
\multiply\theScratch by 60
\theMinute=\theTime
\advance\theMinute by -\theScratch
\theMinuteTens=\theMinute
\divide\theMinuteTens by 10
\theScratch=\theMinuteTens
\multiply\theScratch by 10
\advance\theMinute by -\theScratch

\def\today{{\number\day\space
 \ifcase\month\or
  January\or February\or March\or April\or May\or June\or
  July\or August\or September\or October\or November\or December\fi
 \space\number\year}}













\define\Bc{{\Cal B}}                         




\define\biggnm#1{                            
  \bigg|\bigg|#1\bigg|\bigg|}











\define\Cpx{\bold C}                         








\define\eh{\hat e}                           

\define\eqdef{{\;\overset\text{def}\to=\;}}     





\define\fpamalg#1{{\dsize\;                  
     \operatornamewithlimits*_{#1}\;}}


\define\freeprodi{\mathchoice                
     {\operatornamewithlimits{\ast}
      _{\iota\in I}}
     {\raise.5ex\hbox{$\dsize\operatornamewithlimits{\ast}
      _{\sssize\iota\in I}$}\,}
     {\text{oops!}}{\text{oops!}}}












\define\Hil{{\mathchoice                     
     {\text{\eusm H}}
     {\text{\eusm H}}
     {\text{\eusms H}}
     {\text{\eusmss H}}}}

























\define\lrnm#1{\left|\left|#1\right|\right|} 

\define\lspan{\text{\rm span}@,@,@,}         



\define\Mksa{{M_k^{\text{\it s.a\.}}}}      

\define\MvN{{\Cal M}}                        




\define\MvNsa{{\MvN_{\text{\it s.a\.}}}}     

\define\nm#1{||#1||}                         


\define\Naturals{{\bold N}}                  







\define\owedge{{                             
     \operatorname{\raise.5ex\hbox{\text{$
     \ssize{\,\bigcirc\llap{$\ssize\wedge\,$}\,}$}}}}}

\define\owedgeo#1{{                          
     \underset{\raise.5ex\hbox
     {\text{$\ssize#1$}}}\to\owedge}}






\define\Pt{{\widetilde P}}                   



\define\QED{\newline                         
            \line{$\hfill$\qed}\enddemo}

\define\Qt{{\widetilde Q}}                   




\define\rank{\text{\rm rank}}                


\define\Real{{\bold R}}                      


\define\Reals{{\bold R}}                     


\define\restrict{\lower .3ex                 
     \hbox{\text{$|$}}}




\define\smd#1#2{\underset{#2}\to{#1}}          

\define\smdb#1#2{\undersetbrace{#2}\to{#1}}    

\define\smdbp#1#2#3{\overset{#3}\to            
     {\smd{#1}{#2}}}

\define\smdbpb#1#2#3{\oversetbrace{#3}\to      
     {\smdb{#1}{#2}}}

\define\smdp#1#2#3{\overset{#3}\to             
     {\smd{#1}{#2}}}

\define\smdpb#1#2#3{\oversetbrace{#3}\to       
     {\smd{#1}{#2}}}

\define\smp#1#2{\overset{#2}\to                
     {#1}}                                     








\define\tocdots                              
  {\leaders\hbox to 1em{\hss.\hss}\hfill}    


\define\Tr{\text{\rm Tr}}                    


\define\Uc{{\Cal U}}                         



\define\Vc{{\Cal V}}                         


\define\volume{{\text{\rm vol}}}             





\define\xh{\hat x}                           











  \newcount\bibno \bibno=0
  \def\newbib#1{\advance\bibno by 1 \edef#1{\number\bibno}}
  \ifnum\mycitestyle=1 \def\cite#1{{\rm[\bf #1\rm]}} \fi
  \def\scite#1#2{{\rm[\bf #1\rm, #2]}}

  \newcount\notasecflag \notasecflag=0
  \def\notasec{\notasecflag=1}

  \newcount\secno \secno=0 \newcount\subsecno
  \def\newsec#1{\procno=0 \subsecno=0 \notasecflag=0
    \advance\secno by 1 \edef#1{\number\secno}
    \edef\currentsec{\number\secno}}
  \def\newsubsec#1{\procno=0 \advance\subsecno by 1 \edef#1{\number\subsecno}
    \edef\currentsec{\number\secno.\number\subsecno}}

  \newcount\appendixno \appendixno=0
  \def\newappendix#1{\procno=0 \notasecflag=0 \advance\appendixno by 1
    \ifnum\appendixno=1 \edef\appendixalpha{\hbox{A}}
      \else \ifnum\appendixno=2 \edef\appendixalpha{\hbox{B}} \fi
      \else \ifnum\appendixno=3 \edef\appendixalpha{\hbox{C}} \fi
      \else \ifnum\appendixno=4 \edef\appendixalpha{\hbox{D}} \fi
      \else \ifnum\appendixno=5 \edef\appendixalpha{\hbox{E}} \fi
      \else \ifnum\appendixno=6 \edef\appendixalpha{\hbox{F}} \fi
    \fi
    \edef#1{\appendixalpha}
    \edef\currentsec{\appendixalpha}}

  \newcount\procno \procno=0
  \def\newproc#1{\advance\procno by 1
   \ifnum\notasecflag=0 \edef#1{\currentsec.\number\procno}
   \else \edef#1{\number\procno}
   \fi}

  \newcount\tagno \tagno=0
  \def\newtag#1{\advance\tagno by 1 \edef#1{\number\tagno}}
